\documentclass[12pt]{article}
\usepackage{amssymb}
\usepackage{amsmath}
\usepackage{graphicx}
\usepackage{hyperref}
\usepackage{indentfirst}
\usepackage{combelow}
\usepackage{enumerate}
\newcommand{\R}{\mathbb{R}}
\newcommand{\N}{\mathbb{N}}
\def\argmin{\mathop{\rm argmin}\limits}  
\def\argmax{\mathop{\rm argmax}}
\newcommand{\ep}{\hfill $\square$}
\newcommand{\cl}{\mathrm{cl}}
\newcommand{\conv}{\mathrm{conv}}
\renewcommand{\exp}{\mathrm{exp}}
\newcommand{\ext}{\mathrm{ext}}
\newcommand{\bd}{\mathrm{bd}}

    \makeatletter
\def\@fnsymbol#1{\ensuremath{\ifcase#1\or \dagger\or \ddagger\or
   \mathsection\or \mathparagraph\or \|\or **\or \dagger\dagger
   \or \ddagger\ddagger \else\@ctrerr\fi}}
    \makeatother

\newtheorem{theorem}{Theorem}[section]
\newtheorem{definition}{Definition}[section]
\newtheorem{example}{Example}[section]
\newtheorem{lemma}{Lemma}[section]
\newtheorem{remark}{Remark}[section]
\newtheorem{corollary}{Corollary}[section]
\newtheorem{proposition}{Proposition}[section]
\newcommand{\bp}{\bf Proof. \rm}

\newcommand{\titlename}	
{Reducing the complexity of equilibrium problems and applications to best approximation problems}

\newcommand{\authorname}      {Valerian-Alin Fodor$^1$ \and Nicolae Popovici$^{2,}$\thanks{Meanwhile, professor Nicolae Popovici passed away unexpectedly and prematurely.}}
\newcommand{\pdfauthorname}   {Valerian-Alin Fodor and Nicolae Popovici}

\newcommand{\universityname}  {$^{1,2}$Babe\cb{s}-Bolyai University of Cluj-Napoca, Romania}
\newcommand{\facultyname}     {Faculty of Mathematics and Computer Science}
\newcommand{\departmentname}  {Department of Mathematics}
\newcommand{\pmailaddress}    {valerian.fodor@ubbcluj.ro}

\newcommand{\keywordterms}{...}
\hypersetup{colorlinks = true, linkcolor = blue, anchorcolor =red, citecolor = blue, filecolor = red, urlcolor = blue, pdfauthor=\pdfauthorname, pdftitle=\titlename, pdfkeywords=\keywordterms}
\begin{document}
\title{\LARGE{\titlename}}
\frenchspacing
\newcommand{\institution}{
\universityname\\
\facultyname\\
\departmentname}
\author{\authorname}
\date{
{\footnotesize \institution\\
$\begin{array}{l}
\text{email}^1\text{: \texttt{\href{mailto:\pmailaddress}{\pmailaddress}}}
\end{array}$}}

\maketitle

\begin{abstract}
We consider scalar equilibrium problems governed by a bifunction in a finite-dimensional framework. By using classical arguments in Convex Analysis, we show that under suitable generalized convexity assumptions imposed on the bifunction, the solutions of the equilibrium problem can be characterized by means of extreme or exposed points of the feasible domain. Our results are relevant for different particular instances, such as variational inequalities and optimization problems,  especially for best approximation problems.\\[0.1cm]
\end{abstract}

\section{Introduction and preliminaries}
Throughout this paper  $\R^n$ stands for the $n$-dimensional real Euclidean space,
whose norm $\| \cdot\|$ is induced by the usual inner product $\langle \cdot,\cdot\rangle$.

For all $x, y \in S$, we use the notations 
$$\begin{array}{rcl}
[x,y] &:= & \{(1-t)x+ty \mid t \in [0,1]\,\}, \\    
]x,y[ &:=& \{(1-t)x+ty \mid t \in \,]0,1[\,\}.    
\end{array}$$

Recall that a set $S \subseteq \R^n$ is called convex if $$
[x,y] \subseteq S, \; \forall\, x,y \in S.
$$ 
Of course, this is equivalent to say that 
$$
]x,y[\, \subseteq S, \; \forall\, x,y \in S.
$$

Given a convex set $S \subseteq \R^n$ we denote by 
$$
{\rm ext}\,S = \{x^0 \in S \, \mid \, \forall x,y \in S \;:\;  x^0 = \textstyle{\frac{1}{2}}(x+y)  \;\Rightarrow x=y=x^0 \}
$$
the set of extreme points of $S$. A point $x^0$ is said to be an exposed point of $S$ if there is a supporting hyperplane $H$ which supports $S$ at $x^0$ such that $\{x^0 \} = H \cap S.$ We denote the set of exposed points of $S$ by $$ \exp\,S = \{x^0 \in S \, \mid \, \exists c \in \R^n \text{ such that } \argmin_{x \in S} \langle c,x \rangle = \{x^0 \} \} $$
It is well-known that $\exp\,S \subseteq \ext\,S.$ 
\par 
The convex hull of a set $M \subseteq \R^n$, i.e., the smallest convex set in $\R^n$ containing $M$ is denoted by $\conv M$. 

Next, we recall the following well-known theorems (see for example \cite{Minkowski1911} and \cite{Webster1994}):

\begin{theorem}[Minkowski (Krein-Milman)]\label{kreinmilman}
Every compact convex set in $\R^n$ is the convex hull of its
extreme points.
\end{theorem}

\begin{theorem}[Straszewicz]\label{th:Straszewicz}\index{Straszewicz Theorem}
Every compact convex subset $M$ of $\R^n$ admits the
representation$:$
$$\begin{array}{rcl}
M & = & \cl(\conv (\exp\,M)).
\end{array}$$
\end{theorem}

\begin{definition}\label{de-pba}{\em
Let $S$ be a nonempty subset of $\R^n$ and let $x^* \in \R^n$. A
point $x^0\in S$ is said to be {\it an element of best
approximation to $x^*$ from $S$} (or \emph{a nearest point to
$x^*$ from $S$}) \index{best approximation} \index{nearest point}
if
\begin{eqnarray*}
\|x^0-x^*\| \leq \|x-x^*\|,\; \forall\, x \in S.
\end{eqnarray*}
The problem of best approximation of $x^*$ by elements of $S$ consists in finding all elements of best approximation to $x^*$ from $S$. The solution set 
$$
P_S(x^*) := \{x^0 \in S \mid \; \|x^0-x^*\| \leq \|x-x^*\|,\; \forall\, x \in S\}
$$
is called \emph{the metric projection of $x^*$ on $S$}. 
}\end{definition} 

\begin{remark}{\em
The problem of best approximation is an optimization problem, 
$$\begin{array}{c}
\left\{
\begin{array}{l}
f(x) \longrightarrow \min\\
x \in S,
\end{array}
\right.
\end{array}$$
whose objective function $f : \R^n \to \R$ is defined for all $x \in \R^n$ by
$$
f(x) := \|x-x^*\|.
$$
Actually, we have
$$
P_S(x^*) = \argmin_{x \in S} f(x).
$$
}\end{remark}

\begin{definition}\label{de-fps}{\em
Let $S$ be a nonempty subset of $\R^n$ and let $x^* \in \R^n$, we say that 
 $x^0\in S$ is a farthest point from $S$ to 
$x^*$ if
$$
\|x^0-x^*\| \geq \|x-x^*\|,\; \forall\, x \in S, $$
i.e.,
$$ x^0 \in \argmax_{x \in S} \|x-x^*\|.
$$
}\end{definition}

In this paper we will use the following well known results from Convex Analysis (see for example \cite{BrecknerPopovici2006}).

\begin{proposition}
Any farthest point from a nonempty set $S \subseteq \R^n$ to a point $x^* \in \R^n$ is an exposed point of $S,$ i.e., 
$$\argmax_{x \in S}  \|x-x^*\| \subseteq \exp\, S.$$

\end{proposition}

\begin{theorem}[existence of elements of best approximation]\label{existentapcmba}\,\\
If $S$ is a nonempty closed subset of $\R^n$, then for every $x^*\in\R^n$ there is an element of best
approximation to $x^*$ from $S$. In other words, we have
$$
P_S(x^*) \neq \emptyset,\; i.e., \; {\rm card} (P_S(x^*)) \geq 1.
$$
\end{theorem}

\begin{theorem}[unicity of the element of best approximation]\label{unicitatepcmba}\,\\
If $S \subseteq \R^n$ is a nonempty convex set and $x^*\in\R^n$,
then there exists at most one element of best approximation to
$x^*$ from $S$. In other words, we have
$$
{\rm card} (P_S(x^*)) \leq 1.
$$
\end{theorem}

\begin{theorem}[characterization of elements of best approximation]\label{sufneccmba}
Let $S\subseteq\R^n$, let $x^0\in S$, and let $x^*\in\R^n$. Then
the following hold:
\begin{itemize}
\setlength{\itemsep}{-0.5pt}
\item[{\rm (a)}] If $\langle x-x^0,x^*-x^0\rangle \leq 0$ for all
$x\in S$, then $x^0$ is an element of best approximation to $x^*$ from $S$.
\item[{\rm (b)}] If $S$ is convex and $x^0$ is an element
of best approximation to $x^*$ from $S$, then we have that
$\langle x-x^0,x^*-x^0\rangle \leq 0$ for all $x\in
S$.
\end{itemize}
\end{theorem}

\begin{corollary}\label{co-charP_S(x*)}
Let $S\subseteq\R^n$ be a nonempty convex set and let $x^*\in\R^n$. Then

$$
P_S(x^*) = \{ x^0 \in S \mid \langle x-x^0,x^*-x^0\rangle \leq 0,\; \forall\, x\in
S\}.
$$
\end{corollary}

This paper is organized as follows. In Section \ref{bestapp}, we proved that if $S$ is a nonempty convex set, then the elements of best approximation to an arbitrary element in $\R^n$ from $S$, can be characterized by means of the Gauss map (Remark \ref{GaussMap}). If $S$ is also closed, then we have obtained that $x^0$ is the element of best approximation to an arbitrary element $x^*$ from $S$, if and only if $x^*$ is an element of the translated normal cone to $S$ at $x^0$ by the vector $x^0$ (Proposition \ref{charofmetricprojbynormalcones}). This fact led us to Proposition \ref{partitionofR^nbynormalcones}, 
where we have proved that the set $\{x + N_S (x) \setminus \{0\} \mid x \in \bd \, S \}$ is a partition of $\R^n \setminus S$.

In Section \ref{Equilibriumproblems}, we consider the bifunction $g: A \times \conv M \to \R$ which is quasiconvex in the second argument. For such a bifunction, we show that the equilibrium points of g are precisely the equilibrium points of 
$g|_{A \times M}$. We also point out the particular case $M=\ext S$ for some Minkowski set $M$. 
Furthermore, Theorem \ref{th-eq-h-qc} also led us to our main result, Theorem \ref{chareqbyexposed}, where we show that if $g: S \times S$ is quasiconvex and lower semicontinous in the second argument then, the equilibrium points of $g$ are precisely the equilibrium points of $g|_{S \times \exp S}$, whenever $S$ is a nonempty compact convex set.

\section{The inverse images of the metric projection}\label{bestapp}

\begin{remark}\label{remarknormalcone}{\em
From a geometric point of view, the property 
$\langle x-x^0,x^*-x^0\rangle \leq 0$ for all $x\in
S$ in assertion (b) of Theorem \ref{sufneccmba} shows that  
$x^*-x^0$ belongs to the so-called normal cone to $S$ at $x^0$, i.e., 
$$N_S(x^0) = \{d \in \R^n \mid  \langle x-x^0,d \rangle \leq 0, \; \forall\, x\in
S\}.$$
}\end{remark}

For $S \subseteq \R^n$ a nonempty convex set, we recall  the \textit{Gauss map} of $S$, introduced in \cite{MartinezPintea2021}, which is a set-valued map, defined as follows:
$$G_S:\R^n\rightrightarrows S^{n-1},\;\;\;\;\; G_S(x) := N_S(x) \cap S^{n-1}.$$

\begin{remark}
\label{GaussMap}
{\em
Let $S \subseteq \R^n$ be a nonempty convex set, let $x^* \in \R^n \setminus S$. Then $x^0$ is an element of best approximation to $x^*$ from $S$ if and only if $$\frac{x^*-x^0}{\|x^*-x^0\|} \in G_S(x^0).$$ 
\\Indeed, $x^0$ is an element of best approximation to $x^*$ from $S$ if and only if 
$$
\begin{array}{rcl}
\langle x-x^0,x^*-x^0\rangle \leq 0, \, \forall \, x \in S & \Leftrightarrow & \left\langle x-x^0,\cfrac{x^*-x^0}{\|x^*-x^0\|} \right\rangle \leq 0, \, \forall \, x \in S\\
& \Leftrightarrow & \cfrac{x^*-x^0}{\|x^*-x^0\|} \in N_S(x^0)\\
& \Leftrightarrow & \cfrac{x^*-x^0}{\|x^*-x^0\|} \in G_S(x^0).
\end{array}
$$
}\end{remark}
Remark \ref{GaussMap} combined with the results of \cite{MartinezPintea2021} has some potential for further developments which are left for a forthcoming paper.

Let $S \subseteq \R^n$ be a nonempty closed convex set. By Theorems \ref{existentapcmba} and \ref{unicitatepcmba} it follows that, for all $x^* \in \R^n$, $P_S(x^*)$ is a singleton. So, in this case, $P_S$ can be considered as a single valued mapping. 

\begin{proposition}\label{charofmetricprojbynormalcones}{\em
Let $S \subseteq \R^n$ be a nonempty closed convex set. Then, for all $x^* \in \R^n$, we have that $P_S(x^*)=x^0$ if and only if 
$x^* \in x^0 + N_S(x^0).$
}\end{proposition}
\bp
Let $x^* \in \R^n .$ Since $S$ is closed and convex, there exists $x^0 \in S$ such that $P_S(x^*)=x^0$. It follows by Corollary \ref{co-charP_S(x*)} that $$\langle x-x^0,x^*-x^0\rangle \leq 0,\; \forall\, x\in
S $$ which, by Remark \ref{remarknormalcone}, is equivalent to $$ x^*-x^0 \in N_S (x^0) \iff x^* \in x^0 + N_S (x^0).$$
\ep

\begin{proposition}\label{partitionofR^nbynormalcones}{\em
Let $S \subseteq \R^n$ be a nonempty closed convex set. Then, the set 
\begin{equation}\label{eq05.08.2022.1}
 \{x + N_S (x) \setminus \{0\} \mid x \in \bd \, S \}   
\end{equation}
is a partition of $\R^n \setminus S$.
In other words, $$\R^n \setminus S =\bigcup\limits_{x \in \bd \, S} (x + N_S (x) \setminus \{0\})$$ and $(x + N_S (x) \setminus \{0\}) \cap (y + N_S (y) \setminus \{0\}) \neq \emptyset$  implies $x=y$.
}\end{proposition}
\bp
Let $x^* \in \R^n \setminus S$ and $x^0 \in S$ such that $P_S(x^*)=x^0$. By Proposition   \ref{charofmetricprojbynormalcones}, we obtain that  $x^* \in x^0 + N_S (x^0) \setminus \{0\}$, yet $$ x^0 + N_S (x^0) \setminus \{0\} \subseteq \bigcup\limits_{x \in S} (x + N_S (x) \setminus \{0\}).$$ 
Subsequently, we get that $$\R^n \setminus S \subseteq \bigcup\limits_{x \in S} (x + N_S (x) \setminus \{0\}).$$

In order to prove the converse, let $x \in S$ and $u \in N_S (x) \setminus \{0\}$. If we suppose that $x+u \in S$ then, by the definition of $N_S (x) \setminus \{0\}$, $$\langle x+u-x,u\rangle = \langle u,u\rangle \leq 0 \implies u=0 $$ which contradicts the fact that $u \in N_S (x) \setminus \{0\}$. Thus $x +N_S (x) \setminus \{0\} \subseteq \R^n \setminus S$, for all $x \in S$, i.e., $$\bigcup\limits_{x \in S} (x + N_S (x) \setminus \{0\}) \subseteq \R^n \setminus S.$$

Therefore, we have proved that $$\R^n \setminus S =\bigcup\limits_{x \in  S} (x + N_S (x) \setminus \{0\}).$$ However, since $x + N_S (x) \setminus \{0\}$ is nonempty if and only if $x$ is an boundary point of $S$, we obtain that $$\R^n \setminus S =\bigcup\limits_{x \in \bd \, S} (x + N_S (x) \setminus \{0\}).$$

If $(x + N_S (x) \setminus \{0\}) \cap (y + N_S (y) \setminus \{0\}) \neq \emptyset$, then there is an $u \in N_S (x) \setminus \{0\}$ and $v \in N_S (y) \setminus \{0\})$ such that
$x+u = y+v$. Furthermore, since $u \in N_S (x) \setminus \{0\}$, we obtain  $\langle y-x,u \rangle = \langle u-v,u \rangle \leq 0$, therefore $\|u\|^2 \leq \langle u,v \rangle$. By similar reasoning, since $v \in N_S (y) \setminus \{0\}$, we obtain $\|v\|^2 \leq \langle u,v \rangle$. Thus, $$0 \leq \|x-y\|^2 = \|u-v\|^2 = \|u\|^2 -2 \langle u,v \rangle + \|v\|^2  \leq 0$$ Henceforth, $\|x-y\|=0 $, which implies $ x = y$.
\ep

\begin{remark}{\em
Alternatively, Proposition \ref{partitionofR^nbynormalcones} can be proved as follows. By Proposition \ref{charofmetricprojbynormalcones}, it follows that for all $x \in S$, the set $x + N_S(x)$ is the inverse image of $x$ through the mapping $P_S$, $P_S^{-1}(x)$. 

If we consider the restriction $P_S|_{\R^n \setminus S}$ of the mapping $P_S$ then, for all $x \in S$, we have that the set $x + N_S (x) \setminus \{0\}$ is the inverse image of $x$ through the mapping $P_S|_{\R^n \setminus S}$.

By the equivalence relation induced by $\mathrm{ker}P_S|_{\R^n \setminus S}$, we obtain that the set $$\{x + N_S (x) \setminus \{0\} \mid x \in \bd \, S \}$$ is a partition of $\R^n \setminus S$.

}
\end{remark}

\begin{figure}[ht]
\begin{center}
		\includegraphics[width=12cm]{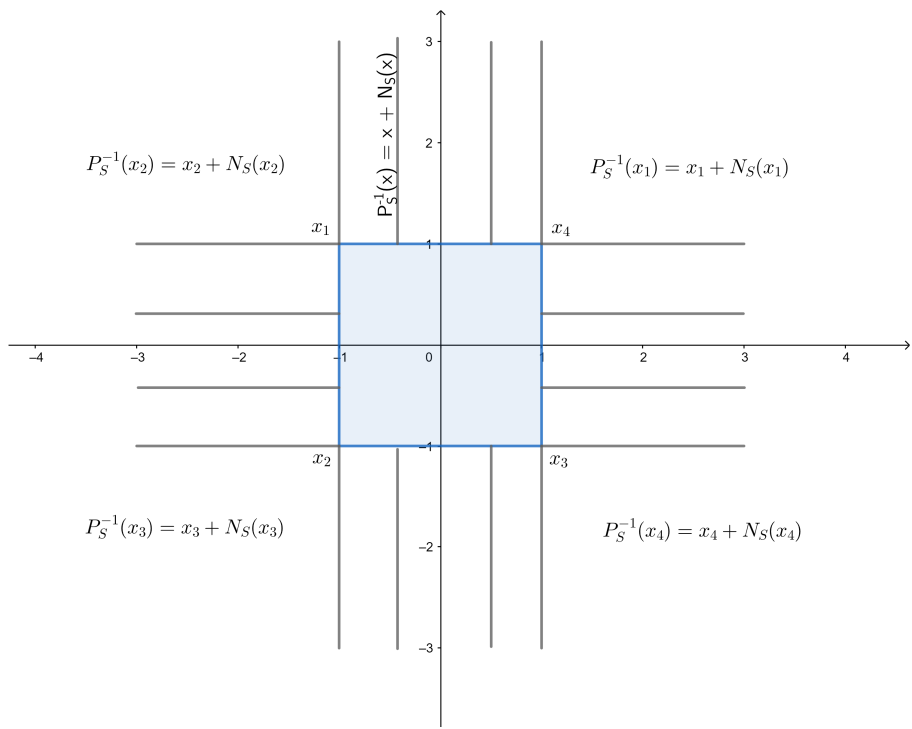}
		\caption{Proposition \ref{partitionofR^nbynormalcones} applied for the  particular case of a square in $\mathbb{R}^2$}
		\label{figure:example22}
\end{center}
\end{figure}

\begin{example}{\em
For $n=2$, consider the set 
$$M = \{x_1=(1,1),x_2=(-1,1),x_3=(-1,-1),x_4=(1,-1)\}\subseteq\R^2$$ and $$S = \conv\,M = [-1,1] \times [-1,1],$$
as in Figure \ref{figure:example22}.
  Obviously, $S$ is a nonempty closed convex subset of $\R^2$ and $$\bd \, S = ]x_1,x_2[ \, \cup \, ]x_2,x_3[ \, \cup \, ]x_3,x_4[ \, \cup \, ]x_1,x_4[  \, \cup \, \{x_1,x_2,x_3,x_4\} .$$
On the four vertices of the square, we have 
$$
\begin{array}{rcl}
N_S (x_1) \setminus \{0\}&=&\{(v_1,v_2) \mid v_1 \geq 0, v_2 \geq 0   \} \setminus \{(0,0)\}\\
N_S (x_2) \setminus \{0\} &=&\{(v_1,v_2) \mid v_1 \leq 0, v_2 \geq 0   \} \setminus \{(0,0)\}\\
N_S (x_3) \setminus \{0\}&=&\{(v_1,v_2) \mid v_1 \leq 0, v_2 \leq 0   \} \setminus \{(0,0)\}\\
N_S (x_4) \setminus \{0\}&=&\{(v_1,v_2) \mid v_1 \geq 0, v_2 \leq 0   \} \setminus \{(0,0)\}.
\end{array}
$$ 
On the other points of the boundary, we have 
$$
N_S (x) \setminus \{0\}=\left\{
\begin{array}{l}
\{(0,v) \mid v > 0 \} \text{, for all }x \in ]x_1,x_2[\\
\{(v,0) \mid v < 0 \} \text{, for all }x \in ]x_2,x_3[\\
\{(0,v) \mid v < 0 \} \text{, for all }x \in ]x_3,x_4[\\
\{(v,0) \mid v > 0 \} \text{, for all }x \in ]x_1,x_4[.
\end{array}
\right.
$$
It is easy to see that  $\R^n \setminus S = \bigcup\limits_{x \in \bd \, S} (x + N_S (x) \setminus \{0\})$.

}
\end{example}

\section{Equilibrium problems}\label{Equilibriumproblems}

The equilibrium problem, introduced in \cite{MuuOettli1992}, has been formulated in a more general way in \cite[p. 18]{KassayRadulescu2019}. We propose a slightly modified definition. Let $g : A \times B \to \R$ be a ``bifunction'', where $A$ and $B$ are nonempty sets. 

\begin{definition}{\em
The \emph{equilibrium problem} with respect to $g : A \times B \to \R$ and a couple of nonempty subsets $A' \subseteq A$ and $B' \subseteq B$, consists in finding the elements
$x^0 \in A'$ satisfying
\begin{eqnarray*}
g(x^0,x) \leq 0 \;\;\mbox{\rm for all}\;\; x \in B'.
\end{eqnarray*}
The set of all solutions of the equilibrium problem 
will be denoted by
$$
\mbox{\rm eq}(g \mid A', B') :=  \{x^0 \in A' \mid g(x^0,x) \leq 0, \; \forall\, x \in B'\}.
$$
}\end{definition}

\begin{remark}{\em 
It is easy to see that 
$$
\mbox{\rm eq}(g \mid A', \emptyset) = A'
$$
and that
$$ 
\mbox{\rm eq}(g \mid A', B') \subseteq \mbox{\rm eq}(g \mid A', B''),\, \forall B'' \subseteq B' 
$$

}\end{remark}

\begin{example}[optimization problems]\label{optprobl} {\em 
Consider a minimization problem 
\begin{eqnarray*}
\left\{
\begin{array}{l}
f(x) \longrightarrow \min\\
x \in S,
\end{array}
\right.
\end{eqnarray*}
where $f : \R^n \to \R$ is a function and $S \subseteq \R^n$ is a nonempty set. 
By defining the bifunction $g : \R^n \times\R^n \to \R$ as
$$
g(u,v) := f(u) - f(v), \; \forall\, (u,v) \in \R^n \times\R^n,
$$
we obtain 
$$
\mbox{\rm eq}(g \mid S, S) = \argmin_{x \in S} f(x).
$$
}\end{example}

\begin{example}[variational inequalities]\label{VI}{\em
Let $T : S \to \R^n$ be a function defined on a nonempty set $S \subseteq \R^n.$ The problem of finding $x^0 \in S$ such that $$ \langle T(x^0), x-x^0 \rangle \geq 0, \forall\, x \in S $$
is called a variational inequality. Denote by sol(VI) the set of its solutions. By defining the bifunction 
$g : S \times S \to \R$ as $$ g(u,v) := \langle T(u),u-v \rangle, \; \forall\, (u,v) \in S \times S,$$ we obtain $$ \mbox{\rm eq}(g \mid S, S) = \mbox{\rm sol(VI)} $$

}\end{example}

\begin{example}[the best approximation problem]\label{bestapprox}{\em $\,$

\begin{enumerate}[a)]
\item The problem of best approximation of $x^*$ by elements of $S$ fits the model described in Example \ref{optprobl}, where $$f(x) = \|x-x^*\| \, \mbox{\rm and} \, g(u,v): = f(u) - f(v)$$ hence 
$$\mbox{\rm eq}(g \mid S, S) = \argmin_{x \in S} \|x-x^* \|. $$

\item Another way of seeing the best approximation problem as an equilibrium problem is to consider the bifunction $g :\R^n \times \R^n \to \R$ defined  by
$$
g(u,v) := \langle v-u,x^* - u \rangle, \; \forall\, (u,v) \in \R^n \times \R^n.
$$
According to Theorem \ref{sufneccmba}, 
$$
\mbox{\rm eq}(g \mid S, S) = \{x^0 \in S \mid \langle x-x^0,x^*-x^0\rangle \leq 0,\; \forall\, x \in S\} \subseteq P_S(x^*),
$$
the equality being true whenever $S$ is convex, i.e., $$ \mbox{\rm eq}(g \mid S, S) = P_S(x^*).$$
Actually, by considering the function $T : S \to \R^n$ defined by $T(x) = x-x^*$ for all $x \in S,$ we recover
$$g(u,v) = \langle T(u),u-v \rangle, \forall \, (u,v) \in S \times S, $$ 
hence, under the convexity assumption on $S$ we can reduce the best approximation problem to a variational inequality:
$$P_S(x^*) = \mbox{\rm sol(VI)}. $$
\end{enumerate}
}\end{example}

\begin{example}[the farthest point problem]\label{farpr}{\em  Let $S \subseteq \R^n$ be a nonempty set and let $x^* \in \R^n.$ The problem of finding the farthest points from $S$ to $x^*$ fits the model described in Example \ref{optprobl}, where $$f(x) = - \|x-x^*\| \, \mbox{\rm and} \, g(u,v): = f(u)-f(v),$$  hence
$$\mbox{\rm eq}(g \mid S, S) = \argmax_{x \in S} \|x-x^*\|$$

}\end{example}

\begin{theorem}\label{th-eq-h-qc}
Let $A$ be a nonempty set, let $S = \conv M$ for some nonempty set $M \subseteq \R^n$ and let $g: A \times S \to R$ be a bifunction. If for every $u \in A$, the function $h = g(u,\cdot) : S \to \R$ is quasiconvex,
i.e.,
$$
h((1-t) v' + t v'') \leq \max\{h(v'),h(v'')\}
$$
for all $v,v' \in S$ and $t \in [0,1]$, then
$$
\mbox{\rm eq}(g \mid A, S) =  \{x^0 \in A \mid g(x^0,x) \leq 0, \; \forall\, x \in M\}, 
$$ i.e.,
$$ \mbox{\rm eq}(g \mid A, S)=\mbox{\rm eq}(g \mid A, M).
$$
\end{theorem}
\bp
We denote by 
$$
E: = \{x^0 \in A \mid g(x^0,x) \leq 0, \; \forall\, x \in M\}.
$$
It is obvious that $\mbox{\rm eq}(g \mid A, S) \subseteq E$.
In order to prove the converse, let $x^0 \in E$ and $x \in S$ arbitrary chosen. 
Since $x \in S = \conv\,M $, this implies that there exists  $ k \in \mathbb{N}$, $ x^1,x^2,...,x^k \in M$ and $t_1,t_2,...,t_k \geq 0 $ such that $\sum_{i=1}^{k} t_i = 1$ and that $x = \sum_{i=1}^{k} t_i x^i$. 
Therefore, given that $g(x^0, \cdot)$ is quasiconvex, we have 
$$
g(x^0,x) = g(x^0,\sum_{i=1}^{k} t_i x^i) \leq \max \{ g(x^0,x^i) \mid i =1,\dots,k \} \leq 0,
$$
since $x^1,x^2,...,x^k \in M$.
\ep

\begin{remark}\label{quasiconvexityreducestoquasiconcavity}{\em
Consider the minimization problem described in Example \ref{optprobl}, where $g(u,v) = f(u) - f(v) \,
\mbox{\rm for all} \, u,v \in S.$ Since for every $u \in S$ we have $\; g(u, \cdot) = f(u) - f,$  the quasiconvexity of $g(u, \cdot)$ for some $u \in S$ reduces to the the quasiconcavity of $f.$ Thus we deduce from the Theorem \ref{th-eq-h-qc} the following result.
}\end{remark}
\begin{corollary}\label{minquasiconcave}
Assume that $S = \conv\, M$ for some nonempty set $M \subseteq \R^n.$ If $f : S \to \R$ is a quasiconcave function, then $$\argmin_{x \in S} f(x) = \{ x^0 \in S \mid f(x^0) \leq f(x), \forall \, x \in M  \} \supseteq \argmin_{x \in M} f(x).$$ 
Moreover, $\argmin_{x \in S} f(x)$ is nonempty if and only if so is $\argmin_{x \in M} f(x),$ hence 
$$\min\, f(S) = \min\, f(M).$$ 
\end{corollary}

\begin{remark}{\em 
The assumptions on quasiconvexity of $g(u,\cdot)$ in Theorem \ref{th-eq-h-qc} and quasiconcavity of $f$ in Corollary \ref{minquasiconcave} are essential, as shown by the next example.
}\end{remark}

\begin{example}{\em

Let $n=1$, $ M=\{-1,1\}$, $ S = \conv\, M=[-1,1]$.  Consider the function 
$$
\left\{\begin{array}{l}
f:S \to \R\\
f(x)=x^2
\end{array}\right.$$
 and the bifunction 
 $$
\left\{\begin{array}{l}
g:  \R^2 \to \R\\
g(u,v) =f(u)-f(v), \; \forall\, (u,v) \in \R^2.
\end{array}\right.
$$
Clearly, the $f$ is not quasiconcave, hence function $ g(u,\cdot) : S \to \R$  is not quasiconvex for any $ u \in S$, in view of Remark \ref{quasiconvexityreducestoquasiconcavity}. It is easy to see that 
$$
\begin{array}{rcl}
\mbox{\rm eq}(g \mid S, S) &=& \argmin_{x \in S} f(x)\\
& = & \{0\} \\
&\nsupseteq&\mbox{\rm eq}(g \mid S, M)\\
& =& \{ x^0 \in S  \mid g(x^0,x) \leq 0, \forall \; x \in M \}\\
& =& S.
\end{array}
$$
Of course, this example also shows that the quasiconcavity assumption imposed on $f$ in Corollary \ref{quasiconvexityreducestoquasiconcavity} is essential, because 
$$
\argmin_{x \in S} f(x) = \{0\} \nsupseteq \argmin_{x \in M} f(x) = M.
$$
}
\end{example}


\begin{remark}{\em
Under the hypothesis of Corollary \ref{minquasiconcave}, the inclusion $$\argmin_{x \in S} f(x) \subseteq \argmin_{x \in M} f(x)$$ does not hold in general, as shown by the following example.
}\end{remark}

\begin{example}{\em

Let $n=1$, $ M=\{-1,1\}$, $ S = \conv M=[-1,1]$.  Consider the function $f : S \to \R$ defined as $f(x) = \mbox{\rm max} \{0,x\},$ for all $x \in S.$ Obviously, f is nondecreasing, hence quasiconcave. However, 
$$\argmin_{x \in S} f(x) = [-1,0] \nsubseteq \argmin_{x \in M} f(x) = \{-1\}.$$

}
\end{example}

\begin{corollary}\label{solvarineq} 
Assume that $S = \conv\, M$ for some nonempty set $M \subseteq \R^n.$ and let $T : S \to \R^n$ be an arbitrary function. Then the set of solutions $$\mbox{\rm sol(VI)} : = \{ x^0 \in S \mid \langle T(x^0), x-x^0 \rangle \geq 0, \forall \, x \in S \}$$ to the variational inequality introduced in Example \ref{VI}, admits the following representation $$\mbox{\rm sol(VI)} : = \{ x^0 \in S \mid \langle T(x^0), x-x^0 \rangle \geq 0, \forall \, x \in M \}.$$
  
\end{corollary}

In the book by Breckner and Popovici \cite[C 5.2.7, p. 82]{BrecknerPopovici2006} we have the following remark:

\begin{remark}[Minkowski]{\em
 Let $S \subseteq \R^n$ be a compact convex set. Then, for each
subset $M$ of $S$, the following equivalence holds$:$
\begin{eqnarray*}
S = \conv\, M & \Longleftrightarrow & {\rm ext}\,S \subseteq M.
\end{eqnarray*}
}\end{remark}

\begin{corollary}\label{cor:Krein}
If $S \subseteq \R^n$ is a nonempty convex compact set and function $g(u,\cdot)$ is quasiconvex on $S$ for every $u \in A$, then
$$
\mbox{\rm eq}(g \mid A, S) =  \{x^0 \in A \mid g(x^0,x) \leq 0, \; \forall\, x \in {\rm ext} S\}
$$
i.e.,
$$
\mbox{\rm eq}(g \mid A, S) = \mbox{\rm eq}(g \mid A, {\rm ext} S).
$$
\end{corollary}
\bp
Follows by Theorem \ref{th-eq-h-qc} and Minkowski Theorem  \ref{kreinmilman}.
\ep

\begin{remark}{\em 
The conclusion of Corollary \ref{cor:Krein} still holds if $S$ is a so-called ``Minkowski set'' (sets which were introduced in \cite{MartinezPintea2016}, i.e., closed, possibly unbounded, sets which can be represented as the convex hull of their extreme points).
}\end{remark}

\begin{remark}{\em 
The hypothesis of closeness in Corollary \ref{cor:Krein} is crucial as it shown in the next example.
}\end{remark}

\begin{example}{\em
Let $n=1, S =]-1,1[.$  Consider the function $f:S \to \R$ defined by  $f(x)=-x^2$ and the bifunction $g:  \R^2 \to \R$ defined by 
$$
g(u,v) :=f(u)-f(v), \; \forall\, (u,v) \in \R^2.
$$
Clearly, $\forall \; u \in S$ the function $ g(u,\cdot) : S \to \R$  is quasiconvex (even convex). It is easy to see that 
$$
\mbox{\rm eq}(g \mid S, S) =  \argmin_{x \in S} f(x) = \emptyset \neq \mbox{\rm eq}(g \mid S, {\rm ext} S) = \mbox{\rm eq}(g \mid S, \emptyset) = S. 
$$
}\end{example}


\begin{theorem}\label{chareqbyexposed}
If $S \subseteq \R^n$ is a nonempty convex compact set and the function $g(u,\cdot)$ is quasiconvex and lower semicontinous on $S$ for every $u \in S$, then
$$
\mbox{\rm eq}(g \mid S, S) =  \{x^0 \in S \mid g(x^0,x) \leq 0, \, \forall\, x \in  \exp\, S\}=\mbox{\rm eq}(g \mid S,\exp\, S).
$$
\end{theorem}
\bp
Let $x^0 \in S$ such that $g(x^0,x) \leq 0, \forall \, x \in \exp\,S.$ 
We prove that 
$$
g(x^0,y) \leq 0, \forall \, y \in S.
$$
Let $y \in S.$ By Theorem \ref{th-eq-h-qc}, it follows that
$$
\begin{array}{ll}
&\{x^0 \in S \mid g(x^0,x) \leq 0, \; \forall\, x \in \exp\,S\}\\
 \subseteq &\{x^0 \in S \mid g(x^0,x) \leq 0, \; \forall\, x \in \conv(\exp\, S)\}.
\end{array} 
$$
Hence 
\begin{equation}\label{eq:steluta}
g(x^0,x) \leq 0, \; \forall\, x \in \conv(\exp\,S).
\tag{$\ast$}\end{equation}
By Straszewicz's theorem (Theorem \ref{th:Straszewicz}), $S =  \cl(\conv(\exp\,S))$ and, since $y \in S,$ there exists a sequence $(y^k)_{k \in \N}$ in $\conv(\exp\,S)$ which converges to $y.$\\
According to \eqref{eq:steluta}, we have $g(x^0,y^k) \leq 0, \forall \,  k \in \N$,  i.e., $$y^k \in L: = \{  z \in S  \mid g(x^0,z) \leq 0   \} ,\; \forall \, k \in \N.$$
Since $g(x^0,\cdot) : S \to \R$ is lower semicontinuous and $S$ is closed, it follows that the level set $L$ is closed with respect to the induced topology in $S$ from $\R^n$, hence 
$$
y = \lim_{k\to\infty} y^k \in \cl\,L = L.
$$
Thereby, $g(x^0,y) \leq 0$ and, since $y$ was arbitrary chosen from $ S$, we obtain that $$\mbox{\rm eq}(g \mid S, S) \supseteq  \{x^0 \in S \mid g(x^0,x) \leq 0, \; \forall\, x \in \exp\,S\}.$$ The reverse inclusion is obvious.

\ep

\begin{corollary}\label{cor:minquasiconcave2}
Assume that $S = \mbox{\rm cl}({\conv\, M})$ for some nonempty set $M \subseteq \R^n.$ If $f : S \to \R$ is a quasiconcave upper semicontinous function, then $$\argmin_{x \in S} f(x) = \{ x^0 \in S \mid f(x^0) \leq f(x), \forall \, x \in M  \} \supseteq \argmin_{x \in M} f(x).$$ 
Moreover, $\argmin_{x \in S} f(x)$ is nonempty if and only if so is $\argmin_{x \in M} f(x),$ hence  
$$\min\, f(S) = \min\, f(M)$$ 
\end{corollary}



\end{document}